\documentclass[a4paper,12pt]{article}
\usepackage{amscd}
\usepackage{amsmath,amsfonts,amssymb,amscd}
\usepackage{indentfirst,graphicx,epsfig}
\usepackage{graphicx,psfrag}
\input{epsf}

\newtheorem{thm}{Theorem}

\newtheorem{lem}{Lemma}

\newtheorem{remark}{Remark}
\newtheorem{cor}{Corollary}

\def\qed{\hfill \nopagebreak\rule{5pt}{8pt}}
\def\pf{\noindent {\it Proof.} }

\title{\bf\large The generalized $3$-connectivity of random
graphs\footnote{Supported by NSFC and the ``973" program. } }
\author{
\small  Ran Gu, Xueliang Li, Yongtang Shi\\
\small Center for Combinatorics and LPMC-TJKLC \\
\small Nankai University, Tianjin 300071, China \\
\small guran323@163.com, lxl@nankai.edu.cn,  shi@nankai.edu.cn
\date{}}
\begin{document}
\maketitle
\begin{abstract}
The generalized connectivity of a graph $G$ was introduced by
Chartrand et al. Let $S$ be a nonempty set of vertices of $G$, and
$\kappa(S)$ be defined as the largest number of internally disjoint
trees $T_1, T_2, \cdots , T_k$ connecting $S$ in $G$. Then for an
integer $r$ with $2 \leq r \leq n$, the {\it generalized
$r$-connectivity} $\kappa_r(G)$ of $G$ is the minimum $\kappa(S)$
where $S$ runs over all the $r$-subsets of the vertex set of $G$.
Obviously, $\kappa_2(G)=\kappa(G)$, is the vertex connectivity of
$G$, and hence the generalized connectivity is a natural
generalization of the vertex connectivity. Similarly, let
$\lambda(S)$ denote the largest number $k$ of pairwise edge-disjoint
trees $T_1, T_2, \ldots , T_k$ connecting $S$ in $G$. Then the {\it
generalized $r$-edge-connectivity} $\lambda_r(G)$ of $G$ is defined
as the minimum $\lambda(S)$ where $S$ runs over all the $r$-subsets
of the vertex set of $G$. Obviously, $\lambda_2(G) = \lambda(G)$.

In this paper, we study the generalized $3$-connectivity of random
graphs and prove that for every fixed integer $k\geq 1$,
$$p=\frac{{\log n+(k+1)\log \log n -\log \log \log n}}{n}$$ is a sharp
threshold function for the property $\kappa_3(G(n, p)) \geq k$,
which could be seen as a counterpart of Bollob\'{a}s and Thomason's
result for vertex connectivity. Moreover, we obtain that $\delta
(G(n,p)) - 1 = \lambda (G(n,p)) - 1 = \kappa (G(n,p)) - 1 \le
{\kappa _3}(G(n,p)) \le {\lambda _3}(G(n,p)) \le \kappa (G(n,p)) =
\lambda (G(n,p)) = \delta (G(n,p))$ almost surely holds, which
could be seen as a counterpart of Ivchenko's result.\\[2mm]
[2mm] AMS Subject Classification (2010): 05C40, 05C80, 05D40.
\end{abstract}

\section{Introduction}

All graphs considered here are finite, undirected, and have no loops
or multiple edges. For standard graph-theoretic notation and
terminology the reader is referred to \cite{Boo}. In particular,
denote by $e[S]$ the number of edges in the induced subgraph by a
set $S$ of vertices. The generalized connectivity of a graph $G$,
introduced by Chartrand et al. in \cite{ec4, COZ}, is a
generalization of the concept of the vertex connectivity. Let $G$ be
a nontrivial connected graph of order $n$ and $r$ an integer with $2
\leq r \leq n$. For a set $S$ of $r$ vertices of $G$, a collection
$T_1, T_2, \cdots, T_k$ of trees in $G$ is said to be internally
disjoint ones connecting $S$ if $E(T_i)\cap E(T_j) = \emptyset$ and
$V(T_i) \cap V (T_j) = S$ for any pair of distinct integers $i$ and
$j$ with $1 \leq i, j\leq k$. Let $\kappa(S)$ denote the largest
number of internally disjoint trees connecting $S$ in $G$. The {\it
generalized $r$-connectivity}, denoted by $\kappa_r(G)$, of $G$ is
then defined by $\kappa_r(G) = \min\{\kappa(S)|~S \subseteq V (G)~
and~ |S| = r\}$. Thus, $\kappa_2(G)=\kappa(G)$, where $\kappa(G)$
denotes the vertex connectivity of $G$. Set $\kappa_r(G)=0$ when $G$
is disconnected. Similarly, one can define the generalized
edge-connectivity. Let $\lambda(S)$ denote the largest number $k$ of
pairwise edge-disjoint trees $T_1, T_2, \ldots, T_k$ in $G$ such
that $S\subseteq V(T_i)$ for every $i$ with $1 \leq i \leq k$. Then
the {\it generalized $r$-edge-connectivity} $\lambda_r(G)$ of $G$ is
defined as $\lambda_r(G)=\min\{\lambda(S)|~S \subseteq V (G)~ and~
|S| = r\}$. Thus, $\lambda_2(G) = \lambda(G)$. Set $\lambda_r(G)=0$
when $G$ is disconnected. The generalized edge-connectivity is
related to an important problem, the Steiner Tree Packing Problem
\cite{Kriesell1, Kriesell2}. There have been many results on the
generalized connectivity and the generalized edge-connectivity, we
refer to the survey \cite{sur} for details.

The generalized connectivity can be motivated by its interesting
interpretation in practice. For example, suppose that $G$ represents
a network. If one wants to connect a set $S$ of vertices of $G$ with
$|S|\geq 3$, then a tree has to be used to connect them. This kind
of tree for connecting a set of vertices is usually called a Steiner
tree, and popularly used in the physical design of VLSI, see
\cite{11}. Usually, one wants to consider how tough a network can
be, for the connection of a set of vertices. Then, the number of
totally independent ways to connect them is a measure for this
purpose.

In this paper, we study the generalized $3$-connectivity of random
graphs. The two most frequently occurring probability models of
random graphs are $G(n,M)$ and $G(n,p)$. The first one consists of
all graphs with $n$ vertices having $M$ edges, in which each graph
have the same probability. The model $G(n,p)$ consists of all graphs
with $n$ vertices in which the edges are chosen independently and
with a same probability $p$.  We say that an event $\mathcal{A}$
happens \textit{almost surely} if its happening probability
approaches $1$ as $n\rightarrow \infty $, i.e.,
$Pr[\mathcal{A}]=1-o_n(1)$. Sometimes, it is addressed as
\textit{a.s.} for short. We will always assume that $n$ is the
variable that tends to infinity.

For a graph property $P$, a function $p(n)$ is called a {\it
threshold function} of $P$ if
\begin{itemize}
\item for every $r(n) = O(p(n))$, $G(n, r(n))$ almost surely satisfies $P$; and

\item for every $r'(n) = o(p(n))$, $G(n, r'(n))$ almost surely
does not satisfy $P$.
\end{itemize}

Furthermore, $p(n)$ is called a {\it sharp threshold function} of
$P$ if there exist two positive constants $c$ and $C$ such that
\begin{itemize}
\item for every $r(n) \geq C\cdot p(n)$, $G(n, r(n))$ almost surely satisfies $P$; and
\item for every $r'(n) \leq c\cdot p(n)$, $G(n, r'(n))$ almost surely
does not satisfy $P$.
\end{itemize}

As well-known, for the vertex connectivity Bollob\'{a}s and Thomason
ontained the following result.
\begin{thm}\cite{Boat}\label{lem3} \
If $k\in\mathbb{N}$ and $y\in \mathbb{R}$ are fixed,  and
$M=\frac{n}{2}(\log n + k\log \log n+y+o(1))\in\mathbb{N}$, then
$$\Pr \left[ {\kappa \left( {G\left( {n,M} \right)} \right) = k}
\right] \to 1 - {e^{ - {e^{ - y/k!}}}}$$

and $$\Pr \left[ {\kappa \left( {G\left( {n,M} \right)} \right) = k
+ 1} \right] \to {e^{ - {e^{ - y/k!}}}}.$$
\end{thm}

Additionally, for the classical connectivity of random graphs,
Ivchenko \cite{Iv} obtained the following result.
\begin{thm}\cite{Iv}\label{lem9} \
If $p\left( n \right) \le \frac{{\log n + k\log \log n}}{n}$ for
some fixed $k$, then \[\Pr \left[ {\kappa (G(n,p(n)) = \lambda
(G(n,p(n)) = \delta (G(n,p(n))} \right] \to 1.\]
\end{thm}

We consider the generalized $3$-connectivity of  random graphs. Our
main result is as follows, which could be seen as a counterpart of
Theorem \ref{lem3}.

\begin{thm}\label{thm1}
Let $k\geq 1$ be a fixed integer. Then $p=\frac{{\log n + (k+1)\log \log
n - \log \log \log n}}{n}$ is a sharp threshold function for the
property $\kappa_3(G(n, p)) \geq k$.
\end{thm}

From Theorem \ref{thm1}, we can obtain the following corollary,
which could be seen as a counterpart of Theorem \ref{lem9}.
\begin{cor}\label{cor1}
Let $p=\frac{{\log n + (k+1)\log \log n - \log \log \log n}}{n}$ for
some fixed $k$. Then, almost surely,
\[k\leq \kappa_3(G(n, p)) \leq \lambda_3(G(n, p))\leq k+1.\]
Moreover,
$\delta (G(n,p))-1 =\lambda(G(n,p))-1=\kappa (G(n,p))-1\le{\kappa _3}(G(n,p))\le
{\lambda _3}(G(n,p))\le\kappa (G(n,p))=\lambda (G(n,p)) =\delta(G(n,p))$ almost surely holds.
\end{cor}

\section{Main results}

Throughout the paper {\it log} always denotes the natural logarithm,
and we assume that $k\geq 1$ is a fixed integer. To establish a
sharp threshold function for a graph property the proof should be
two-fold. We first show one easy direction. The following result is
given by Li et al. in \cite{7}, which will be used later.

\begin{lem}\cite{7}\label{lem1} \
For any connected graph G, $\kappa_3(G)\leq \kappa(G)$. Moreover,
the upper bound is sharp.
\end{lem}

We first prove the following result.
\begin{thm}\label{thm2}
$\kappa_3(G(n, \frac{1}{2}\frac{{\log n + (k+1)\log \log n - \log \log
\log n}}{n})) \leq k-1$ almost surely holds.
\end{thm}

We need the following lemma. We call a property $Q$ \textit{convex}
if $F\subset G\subset H$ and $F$ satisfies $Q$, then $H$ satisfies
$Q$ imply that $G$ satisfies $Q$, where $F,\ G,\ H$ are some graphs.
Set $N=\frac{1}{2}n(n-1)$.
\begin{lem}\cite{Boo}\label{lem2} \
If $Q$ is a convex property and $p(1-p)N\rightarrow \infty $, then
$G(n,p)$ almost surely satisfies $Q$ if and only if for every fixed
$x$, $G(n, M)$ almost surely satisfies $Q$, where $M=\lfloor
pN+x(p(1-p)N)^{1/2}\rfloor$.
\end{lem}

\noindent {\bf Proof of Theorem \ref{thm2}:} Let $p=\frac{{\log n +
(k+1)\log \log n - \log \log \log n}}{n}$ and
$M'=\lfloor\frac{1}{2}pN+x\{\frac{1}{2}p(1-\frac{1}{2}p)N\}^{1/2}\rfloor$
for any $x\in \mathbb{R}$, i.e., $M' = \frac{n}{4}(\log n +
(k+1)\log \log n - \log \log \log n + o(1))$.  It is easy to check
that $\frac{1}{2}p(1-\frac{1}{2}p)N\rightarrow \infty $.

Let $M_1=\frac{n}{2}\left( {\log n + (k-1)\log \log n + y + o\left(
1 \right)} \right)\in\mathbb{N}$. By Theorem \ref{lem3}, we have
$$\Pr \left[ {\kappa \left( {G\left( {n,{M_1}} \right)} \right) = k
- 1} \right] \to 1 - {e^{ - {e^{ - y/\left( {k - 1} \right)!}}}}.$$
Hence, for any $\varepsilon >0$, there exists an $N'\in\mathbb{N}$
and a $Y\in\mathbb{R^+}$, such that for any $y<-Y$, $$1 - {e^{ -
{e^{ - y/\left( {k - 1} \right)!}}}} - \Pr \left[ {\kappa \left(
{G\left( {n,{M_1}} \right)} \right) = k - 1} \right] <
\frac{\varepsilon }{2}\ \  \mathrm{and} \ \ {e^{ - {e^{ - y/\left(
{k - 1} \right)!}}}} < \frac{\varepsilon }{2}\,.$$

On the other hand,  there exists an integer $N_1\in\mathbb{N}$, such
that for any $n>N_1$, $M'<M_1$. We have
\begin{eqnarray*}
\Pr \left[ {\kappa \left( {G\left( {n,M'} \right)} \right) \le k - 1} \right]&= &
\sum\limits_{i = 0}^\infty  {\Pr \left[ {\kappa \left( {G\left( {n,M'} \right)} \right)
\le k - 1|\kappa \left( {G\left( {n,{M_1}} \right)} \right) = i} \right]}\\
\null &\null& \cdot \Pr \left[ {\kappa \left( {G\left( {n,{M_1}} \right)} \right) = i} \right]\\
\null  &\geq & \Pr \left[ {\kappa \left( {G\left( {n,M'} \right)} \right)
\le k - 1|\kappa \left( {G\left( {n,{M_1}} \right)} \right) = k - 1} \right]\\
\null  &\null& \cdot \Pr \left[ {\kappa \left( {G\left( {n,{M_1}} \right)} \right) = k - 1} \right]\\
\null  &=&  \Pr \left[ {\kappa \left( {G\left( {n,M'} \right)} \right)
\le k - 1,\kappa \left( {G\left( {n,{M_1}} \right)} \right) = k - 1} \right]\\
\null  &=& \Pr \left[ {\kappa \left( {G\left( {n,{M_1}} \right)}
\right) = k - 1} \right]
\end{eqnarray*}
Hence, we have for any $n>\max \{N',N_1\}$,
\begin{eqnarray*}
1-\Pr \left[ {\kappa \left( {G\left( {n,M'} \right)} \right) \le k - 1} \right]
&\leq& 1-\Pr \left[ {\kappa \left( {G\left( {n,{M_1}} \right)} \right) = k - 1} \right]\\
\null  &<&  {e^{ - {e^{ - y/\left( {k - 1} \right)!}}}} +
\frac{\varepsilon }{2}< \varepsilon.
\end{eqnarray*}
Thus, $\kappa \left( {G\left( {n,M'} \right)} \right) \le k - 1 $
almost surely holds. Obviously, the property that the connectivity
of a given graph is at most $k-1$, is a convex property. By Lemmas
\ref{lem1} and \ref{lem2}, $\kappa_3(G(n, \frac{1}{2}\frac{{\log n +
(k+1)\log \log n - \log \log \log n}}{n})) \leq k-1$ almost surely
holds. \qed

\vspace{3ex}

We leave with the other direction stated below.

\begin{thm}\label{thm3}
$\kappa_3(G(n, \frac{{\log n + (k+1)\log \log n - \log \log \log
n}}{n})) \geq k$ almost surely holds.

\end{thm}

From now on, $p$ is always $\frac{{\log n + (k+1)\log \log n - \log
\log \log n}}{n}$.

\begin{lem}\label{lem7}
$k\leq \kappa(G(n,p))\leq k+1 $ almost surely holds.
\end{lem}

\pf We prove this lemma by a similar method used in the proof of
Theorem \ref{thm2}. Let $M=\lfloor pN+x\{p(1-p)N\}^{1/2}\rfloor$ for
any $x\in \mathbb{R}$, i.e., $M = \frac{n}{2}(\log n + (k+1)\log
\log n - \log \log \log n + o(1))$.  It is easy to check that
$p(1-p)N\rightarrow \infty $. Let $M_0=\frac{n}{2}\left( {\log n +
(k+1)\log \log n + y + o\left( 1 \right)} \right)\in\mathbb{N}$,
${M_2} = \frac{n}{2}(\log n + k\log \log n+y $ $+ o\left(1
\right))\in\mathbb{N}$. By Theorem \ref{lem3}, \[\Pr \left[ {\kappa
\left( {G\left( {n,{M_2}} \right)} \right) = k } \right] \to 1 -
{e^{ - {e^{ - y/ {k } !}}}}.\] Hence, for any $\varepsilon >0$,
there exists an $N''\in\mathbb{N}$ and a $Y\in\mathbb{R^+}$, such
that for any $y<-Y$, $$1 - {e^{ - {e^{ - y/{k} !}}}} - \Pr \left[
{\kappa \left( {G\left( {n,{M_2}} \right)} \right) = k} \right] <
\frac{\varepsilon }{2}\ \ \mathrm{and} \ \ {e^{ - {e^{ - y/{k} !}}}}
< \frac{\varepsilon }{2}.$$ On the other hand, there exists an
integer $N_2\in\mathbb{N}$, such that for any $n>N_2$, $M>M_2$. We
have
\begin{eqnarray*}
\Pr \left[ {\kappa \left( {G\left( {n,M} \right)} \right) \geq k } \right]&= &
\sum\limits_{i = 0}^\infty  {\Pr \left[ {\kappa \left( {G\left( {n,M} \right)} \right)
\geq k|\kappa \left( {G\left( {n,{M_2}} \right)} \right) = i} \right]}\\
\null &\null& \cdot \Pr \left[ {\kappa \left( {G\left( {n,{M_2}} \right)} \right) = i} \right]\\
\null  &\geq & \Pr \left[ {\kappa \left( {G\left( {n,M} \right)} \right) \geq k |\kappa
\left( {G\left( {n,{M_2}} \right)} \right) = k} \right]\\
\null  &\null& \cdot \Pr \left[ {\kappa \left( {G\left( {n,{M_2}} \right)} \right) = k } \right]\\
\null  &=&  \Pr \left[ {\kappa \left( {G\left( {n,M} \right)} \right) \geq k,\kappa
\left( {G\left( {n,{M_2}} \right)} \right) = k } \right]\\
\null  &=& \Pr \left[ {\kappa \left( {G\left( {n,{M_2}} \right)}
\right) = k } \right]
\end{eqnarray*}
Hence, we get for any $n>\max\{N'',N_2\}$,
\begin{eqnarray*}
1-\Pr \left[ {\kappa \left( {G\left( {n,M} \right)} \right) \geq k } \right]
&\leq& 1-\Pr \left[ {\kappa \left( {G\left( {n,{M_2}} \right)} \right) = k }
\right]\\[3mm]
\null  &<&  {e^{ - {e^{ - y/ {k} !}}}} +
\frac{\varepsilon }{2}< \varepsilon.
\end{eqnarray*}
Thus, $\kappa \left( {G\left( {n,M} \right)} \right) \geq k $
almost surely holds. By Lemma \ref{lem2},  $\kappa(G(n, p)) \geq
k$ almost surely holds.

Similarly, we can prove that $\kappa(G(n,p))\leq k+1 $. By Theorem
\ref{lem3},
\[\Pr \left[ {\kappa \left( {G\left( {n,{M_0}} \right)} \right) = k+1}
\right] \to 1 - {e^{ - {e^{ - y/ {(k+1)} !}}}}.\] Hence, for any
$\varepsilon
>0$, there exists an $N^*\in\mathbb{N}$ and a $Y\in\mathbb{R^+}$, such that for any $y<-Y$,
$$1 - {e^{ - {e^{ - y/ {(k+1)} !}}}} - \Pr \left[ {\kappa \left( {G\left(
{n,{M_0}} \right)} \right) = k+1} \right] < \frac{\varepsilon }{2} \
\ \mathrm{and} \ \ {e^{ - {e^{ - y/ {(k+1)} !}}}} <
\frac{\varepsilon }{2}.$$ On the other hand,  there exists an
integer $N_3\in\mathbb{N}$, such that for any $n>N_3$,
$-\log\log\log n< -Y$. Namely, $M<M_0$, and then
\begin{eqnarray*}
\Pr \left[ {\kappa \left( {G\left( {n,M} \right)} \right) \leq k+1} \right]&= &
\sum\limits_{i = 0}^\infty  {\Pr \left[ {\kappa \left( {G\left( {n,M} \right)} \right)
\leq k+1|\kappa \left( {G\left( {n,{M_0}} \right)} \right) = i} \right]}\\
\null &\null& \cdot \Pr \left[ {\kappa \left( {G\left( {n,{M_0}} \right)} \right) = i} \right]\\
\null  &\geq & \Pr \left[ {\kappa \left( {G\left( {n,M} \right)} \right) \le k+1|\kappa
\left( {G\left( {n,{M_0}} \right)} \right) = k+1} \right]\\
\null  &\null& \cdot \Pr \left[ {\kappa \left( {G\left( {n,{M_0}} \right)} \right) = k+1} \right]\\
\null  &=&  \Pr \left[ {\kappa \left( {G\left( {n,M} \right)} \right) \leq k+1,\kappa
\left( {G\left( {n,{M_0}} \right)} \right) = k+1} \right]\\
\null  &=& \Pr \left[ {\kappa \left( {G\left( {n,{M_0}} \right)}
\right) = k+1} \right].
\end{eqnarray*}
Hence, we have for any $n>\max \{N^*,N_3\}$,
\begin{eqnarray*}
1-\Pr \left[ {\kappa \left( {G\left( {n,M} \right)} \right) \leq k+1} \right] &\leq& 1-\Pr
\left[ {\kappa \left( {G\left( {n,{M_0}} \right)} \right) = k+1}
\right]\\[3mm]
\null  &<&  {e^{ - {e^{ - x/ {(k+1)}!}}}} + \frac{\varepsilon }{2}<
\varepsilon.
\end{eqnarray*}
Thus, $\kappa \left( {G\left( {n,M} \right)} \right) \leq k+1$
almost surely holds. By Lemma 2,  $\kappa(G(n, p)) \leq k+1$ almost
surely holds.\qed\vspace{3ex}

From Lemma \ref{lem7}, we know that the minimum degree of $G(n,p)$
is a.s. at least $k$.

Let $G=G(n,p)$ and $D=\frac{ \log n}{\log\log n}$. Almost surely,
the diameter of $G$ is asymptotically equal to $D$, see for example
\cite{Boo}. We call a vertex $v$ \textit{large} if it is of degree
$d_G(v)\geq \frac{ \log n}{100}$, and \textit{small} otherwise.

\begin{lem}\label{lem4}
Almost surely, there does not exist two small vertices within distance at most
$\frac{3}{4}D$ in $G$.
\end{lem}
\pf Denote by $dist_G(x,y)$ the distance between $x$ and $y$ in $G$.
We have
\begin{eqnarray*}
\null  &\null&   \Pr \left[ {\exists\;x,y \in V(G):~{d_G}(x),~{d_G}(y)
< \frac{{\log n}}{{100}},\;dist_G}(x,y) \le \frac{{3D}}{4} \right]\\
\null &\leq& \left( {\begin{array}{*{20}{c}}
n\\
2
\end{array}} \right)\sum\limits_{j = 1}^{\frac{{3D}}{4}} {\left( {\begin{array}{*{20}{c}}
{n - 2}\\
{j - 1}
\end{array}} \right){p^j}{{\left( {\sum\limits_{i = 0}^{\frac{{\log n}}{{100}}}
{\left( {\begin{array}{*{20}{c}}
{n - \left( {j + 1} \right)}\\
i
\end{array}} \right){p^i}{{\left( {1 - p} \right)}^{n - \left( {i + 1} \right)}}} } \right)}^2}} \\
\null   &\leq&\sum\limits_{j = 1}^{\frac{{3D}}{4}} {n{{\left( {2\log
n} \right)}^j}\{{\sum\limits_{i = 0}^{\frac{{\log n}}{{100}}}
{\left( {\begin{array}{*{20}{c}}
{n - 1}\\
{\frac{{\log n}}{{100}}}
\end{array}} \right)} {p^{\frac{{\log n}}{{100}}}}{{\left( {1 - p} \right)}^{n - 1 -
\frac{{\log n}}{{100}}}}}\}^2 }\\
\null   &\leq& \sum\limits_{j = 1}^{\frac{{3D}}{4}} {n{{\left( {2\log n} \right)}^j}}
\{\frac{{\log n}}{{100}}{\left( {\frac{{ne}}{{\frac{{\log n}}{{100}}}}}
\right)^{\frac{{\log n}}{{100}}}}{{\left( {\frac{{\log n}}{{100}}}
\right)^{\frac{{\log n}}{{100}}}}{e^{ - \frac{{\log n}}{n}
\left( {n - 1 - \frac{{\log n}}{{100}}} \right)}}}\}^2\\
\null  &\leq& \sum\limits_{j = 1}^{\frac{{3D}}{4}} {n{{\left( {2\log n}
\right)}^j}{{\left( {\frac{{\log n}}{{100}}{{\left( {100{e^{1 + o(1)}}}
\right)}^{\frac{{\log n}}{{100}}}} \cdot \frac{1}{n}} \right)}^2}} \\
\null  &\leq& \frac{{3\log n}}{{4\log \log n}}{\left( {2\log n}
\right)^{\frac{{3\log n}}{{4\log \log n}}}}{\left( {\frac{{\log
n}}{{100}}} \right)^2}\frac{1}{n}{n^{\frac{8}{{50}}}}\leq {n^{ -
\frac{9}{{100}}}}.
\end{eqnarray*}
The proof is thus complete.
\qed

\begin{lem}\label{lem5}
For a fixed $t\in\mathbb{N}$ and $0<\alpha<1$, almost surely, there
does not exist a subset $S\subset V(G)$, such that $|S| \leq \alpha
tD$ and $e[S] \geq |S| + t$.
\end{lem}

\pf For convenience, let $s = |S|$. Then we have
\begin{eqnarray*}
\null  &\null& \Pr \left[ {\exists \,S:\;\,s \le \alpha tD,\,\;e[S] \ge s + t} \right] \\
\null &\leq& \sum\limits_{s \le \alpha tD} {\left(
{\begin{array}{*{20}{c}}
n\\
s
\end{array}} \right)} \left( {\begin{array}{*{20}{c}}
{\left( {\begin{array}{*{20}{c}}
s\\
2
\end{array}} \right)}\\
{s + t}
\end{array}} \right){p^{s + t}}  \\
\null  &\leq&   \sum\limits_{s \le \alpha tD} {{{\left(
{\frac{{ne}}{s}} \right)}^s}} {\left( {\frac{{\frac{1}{2}s\left( {s - 1}
\right)e}}{{s + t}}} \right)^{s + t}}{p^{s + t}}\\
\null  &\leq& \sum\limits_{s \le \alpha tD} {{{\left( {\frac{{{e^2}s}}{{2(s + t)}}\log n}
\right)}^s}} {\left( {\frac{{se\log n}}{n}} \right)^t}\\
\null  &\leq& \alpha t\frac{{\log n}}{{\log \log n}}{\left( {{e^{2 + o(1)}}\log n}
\right)^{\alpha t\frac{{\log n}}{{\log \log n}}}}{\left(
{\frac{{e\alpha t\frac{{{{\log }^2}n}}{{\log \log n}}}}{n}} \right)^t}\\
\null  &<& \frac{1}{{{n^{\left( {1 - \alpha  - o(1)} \right)t}}}}.
\end{eqnarray*}
The proof is thus complete. \qed

\begin{remark}\label{rem1}
Let $T$ be a rooted tree of depth at most $3D/4$ and let $v$ be a
vertex not in $T$, but with $b$ neighbors in $T$. Let $S$ consist of
$v$, the neighbors of $v$ in $T$ and the ancestors of these
neighbors. Then $|S| \leq3bD/4 + 1 +1\leq 4bD/5$ and $e[S] = |S| + b
- 2$. It follows from the proof of Lemma \ref{lem5} with
$\alpha=4/5$, $t=16$, that we must have $b\leq 18$, with probability
$1-o(n^{-(1/5-o(1))t})\geq 1-o(n^{-3})$.
\end{remark}

\begin{remark}\label{rem2}
Let $\mathcal{P}$ be a set of at most $k$ vertex disjoint paths and
trees, each containing at most $5D/2$ edges, and let $v$ be a vertex
not in $\mathcal{P}$, but with $c$ neighbors in $\mathcal{P}$. Let
$S=\{v\}\cup V(\mathcal{P})$, $|S|\leq 5kD/2+k+1\leq11kD/4$ and
$e[S]=|S|+c-k-1$. By Lemma \ref{lem5} with $\alpha=1/4$, $t= 11k $,
we deduce that with probability at least $1-o(n^{-3})$, $c\leq
12k+1$.
\end{remark}

We first deal with large vertices. The following lemma points out
that for every pair of large vertices in $V(G)$, there exists a
special subgraph containing them, which can be used to find trees
connecting given vertices. Recall that a \textit{t-ary tree} with a
designated root is a tree whose non-leaf vertices all have exactly
$t$ children.

\begin{lem}\label{lem6}
Let $\varepsilon=\varepsilon(n)=\frac{1}{\log\log n}$. Then, almost
surely, for any pair of large vertices $u$ and $v\in V(G)$, there
exists a subgraph $G_{u,v}$ of $G$ that consists of two vertex
disjoint $\frac{\log n}{101}$-ary trees $T_u$ and $T_v$ rooted at
$u$ and $v$, respectively, each having depth
$(\frac{3}{4}-\varepsilon)D$.
\end{lem}

\pf We will show that for any pair of large vertices $u$ and $v$,
the two trees described in Lemma \ref{lem6} exist with probability
$1-o(n^{-3})$.

Firstly, we grow a tree from $u$ using BFS until it reaches depth
$(\frac{3}{4}-\varepsilon)D$. Then we grow a tree starting from $v$
again using BFS until it reaches depth $(\frac{3}{4}-\varepsilon)D$.

We use the notation ${S_i}^{(x)}$ for the number of vertices at
depth $i$ of the BFS tree rooted at $x$.

As growing $T_u$, when we grow the tree from a vertex $x$ at depth
$i$ to depth $i+1$, there may exist some \textit{bad edges} which
connect $x$ to vertices already in $T_u$.

Remark \ref{rem1} implies that with probability $1-o(n^{-3})$, there
exist at most $18$ bad edges from $x$.

For small vertices, from Lemma 4 we can easily get that in the first
$3D/8$ levels, there exists at most one small vertex at each level
a.s..   Furthermore, once a small vertex appears in the BFS tree,
there will be no small vertex in the subtree rooted at that small
vertex. Though there may be more than one small vertex in depth
$3D/8+1$, the number of them will not exceed the number of branches
at root $u$, since one branch contains at most one small vertex in
depth $3D/8+1$, a.s.. Then in depth $3D/8+2$, the number of small
vertices of that level will be no more than the number of vertices
in the depth $3D/8+1$ contained in the branches which have no small
vertex in the previous levels. For the remaining levels of that BFS
tree, we can conclude the similar result.  And note that there will
exist no small vertex in the following levels of branches which
contain small vertices in depth at least $3D/8+1$ of that BFS tree,
a.s..  Hence the number of small vertices contained in each level is
much smaller than the increase of the number of vertices in each
level. Denote by ${t_i}^{(u)}$ the number of small vertices of depth
$i$. Thus we get the following recursion:
\[{S_{i + 1}}^{(u)} \ge \left( {\frac{{\log n}}{{100}} - 18} \right)
\left( {{S_i}^{(u)} - {t_i}^{(u)}} \right) \ge \frac{{\log
n}}{{101}}{S_i}^{(u)}\] We call the operation of deleting some
vertices from a tree as {\it prune a tree}. It is clear that we can
make the current BFS tree a $\frac{\log n}{101}$-ary tree by
pruning.

Then we grow $T_v$, similarly. The only difference is that now we
also say that an edge is {\it bad} if the other endpoint is in
$T_u$.

Hence, \[{S_{i + 1}}^{(v)} \ge \left( {\frac{{\log n}}{{100}} - 36}
\right)\left( {{S_i}^{(v)} - {t_i}^{(v)}} \right) \ge \frac{{\log
n}}{{101}}{S_i}^{(v)}\]

After pruning, we can obtain the required subgraph $G_{u,v}$. \qed

\vspace{5mm} \noindent{\bf Proof of Theorem \ref{thm3}:} In order to
prove Theorem \ref{thm3}, we will show that for any three vertices,
we can find at least $k$ internally disjoint trees connecting them
in $G$.

Given three vertices $u$, $v$ and $w$, we first assume that they are
all large vertices. With the aid of Lemma \ref{lem6}, construct two
vertex disjoint $\frac{\log n}{101}$-ary trees $T_u$ and $T_v$
rooted at $u$ and $v$, respectively, each having depth
$(\frac{3}{4}-\varepsilon)D$.

For every tree $T$, denote the set of leaves of $T$ by $L(T)$.  Let
$u_1, \ldots, u_{\frac{\log n}{101}}$ ($v_1, \ldots, v_{\frac{\log
n}{101}}$) be the vertices in the first depth of $T_u$ ($T_v$
respectively). For each $u_i$ ($v_i$), denote by
$T_{u_i}$($T_{v_i}$) the subtree of $T_u$ ($T_v$) of depth
$\left({\frac{3}{4} - \varepsilon } \right)D - 1$ rooted at $u_i
$($v_i$), $i=1,\ldots, \frac{\log n}{101}$. Call these $T_{u_i}$
($T_{v_i}$) \textit{vice trees}.

For a fixed $T_{u_i}$, let the random variable $A_i$ denote the
number of edges between $L(T_{u_i})$  and $L(T_v)$. Then $A_i$
follows the binomial distribution, i.e., $A_i\sim Bin\left(
{{{\left( {\frac{{\log n}}{{101}}} \right)}^{\left( {\frac{3}{4} -
\varepsilon } \right)D - 1}} \cdot {{\left( {\frac{{\log n}}{{101}}}
\right)}^{\left( {\frac{3}{4} - \varepsilon } \right)D}},\ p}
\right)$. The expectation value  of $A_i$
\begin{eqnarray*}
\mathbb{E}[A_i]=p{\left( {\frac{{\log n}}{{101}}} \right)^{2\left(
{\frac{3}{4} - \varepsilon } \right)D - 1}} \geq
\frac{{101}}{n}{\left( {\frac{{\log n}}{{101}}} \right)^{2\left(
{\frac{3}{4} - \varepsilon } \right)D}}\geq 101{n^{\frac{1}{2} -
2\varepsilon  - \frac{{9.2\left( {\frac{3}{4} - \varepsilon }
\right)}}{{\log \log
 n}}}}.
\end{eqnarray*}
By Chernoff Bounds,
\begin{eqnarray*}
\Pr \left[ {{A_i} < \frac{{100}}{{101}}\mathbb{E}\left[ {{A_i}} \right]} \right]
&\leq &{e^{ - \frac{1}{2} \times \frac{1}{{{{101}^2}}}\mathbb{E}\left[ {{A_i}} \right]}}
\leq  {e^{ - \frac{1}{2} \times \frac{1}{{{{101}^2}}}\times 101{n^{\frac{1}{2} - 2\varepsilon
- \frac{{9.2\left( {\frac{3}{4} - \varepsilon } \right)}}{{\log \log n}}}}}}\\
\null  &= & {e^{ - \log n\left( {\frac{1}{2} \times \frac{1}{{101}}\frac{1}{{\log n}}{n^{\frac{1}{2}
- 2\varepsilon  - \frac{{9.2\left( {\frac{3}{4} - \varepsilon } \right)}}{{\log \log n}}}}} \right)}}\\
\null  &= &{n^{ - \frac{1}{2} \times \frac{1}{{101}}\frac{1}{{\log
n}}{n^{\frac{1}{2} - 2\varepsilon  - \frac{{9.2\left( {\frac{3}{4} -
\varepsilon } \right)}}{{\log \log n}}}}}}\leq  {n^{ -
{n^{\frac{1}{2} - o\left( 1 \right)}}}}.
\end{eqnarray*}

Now, for a fixed $T_{v_j}$, let $A_{ij}$ denote the number of edges
between $L(T_{u_i})$  and $L(T_{v_j})$. Then $A_{ij}\sim Bin\left(
{{{\left( {\frac{{\log n}}{{101}}} \right)}^{\left( {\frac{3}{4} -
\varepsilon } \right)D - 1}} \cdot {{\left( {\frac{{\log n}}{{101}}}
\right)}^{\left( {\frac{3}{4} - \varepsilon } \right)D - 1}},p}
\right)$. We have
\begin{eqnarray*}
\mathbb{E}[A_{ij}]&=&   {\left( {\frac{{\log n}}{{101}}}
\right)^{2\left( {\frac{3}{4} - \varepsilon } \right)D - 2}} \cdot p
\leq \frac{{{{101}^2}}}{{{{\log }^2}n}}{\left( {\frac{{\log n}}{{101}}}
 \right)^{2\left( {\frac{3}{4} - \varepsilon } \right)D}}\frac{{2\log n}}{n}\\[3mm]
\null  &=&  \frac{{{{101}^2} \times 2}}{{\log n}} \cdot
{n^{\frac{1}{2} - 2\varepsilon - \frac{{9.2\left( {\frac{3}{4} -
\varepsilon } \right)}}{{\log \log n}}}}.
\end{eqnarray*}

Also, we can deduce that $\mathbb{E}[A_{ij}]\geq \frac{{{{101}^2}
}}{{\log n}} \cdot {n^{\frac{1}{2} - 2\varepsilon  -
\frac{{9.2\left( {\frac{3}{4} - \varepsilon } \right)}}{{\log \log
n}}}}$. By applying Chernoff Bounds,
\begin{eqnarray*}
\Pr \left[ {{A_{ij}} > 8\mathbb{E}\left[ {{A_{ij}}} \right]} \right]
&\leq &\frac{{{e^{7\mathbb{E}[{A_{ij}}]}}}}{{{8^{8\mathbb{E}[A_{ij}]}}}}\leq
\frac{{{e^{7\frac{{{{101}^2} \times 2}}{{\log n}} \cdot {n^{\frac{1}{2}
- 2\varepsilon  - \frac{{9.2\left( {\frac{3}{4} - \varepsilon }
\right)}}{{\log \log n}}}}}}}}{{{8^{8\frac{{{{101}^2}}}{{\log n}} \cdot {n^{\frac{1}{2}
- 2\varepsilon  - \frac{{9.2\left( {\frac{3}{4} - \varepsilon } \right)}}{{\log \log n}}}}}}}} \\
\null &= & \frac{{{e^{\frac{{{{101}^2} \times 14}}{{\log n}} \cdot {n^{\frac{1}{2}
- 2\varepsilon  - \frac{{9.2\left( {\frac{3}{4} - \varepsilon }
\right)}}{{\log \log n}}}}}}}}{{{e^{\frac{{{{101}^2} \times 8 \times \log 8}}{{\log n}}
\cdot {n^{\frac{1}{2} - 2\varepsilon  - \frac{{9.2\left( {\frac{3}{4} - \varepsilon }
\right)}}{{\log \log n}}}}}}}}\\
\null &= &{e^{\left( {14 - 8\log 8} \right)\frac{{{{101}^2}}}{{\log n}}\cdot{n^{\frac{1}{2}
- 2\varepsilon  - \frac{{9.2\left( {\frac{3}{4} - \varepsilon } \right)}}{{\log \log n}}}}}}\\
\null &=&   {n^{ - \frac{{2.64 \times {{101}^2}}}{{{{\log }^2}n}}
\cdot {n^{\frac{1}{2} - 2\varepsilon  - \frac{{9.2\left(
{\frac{3}{4} - \varepsilon } \right)}}{{\log \log n}}}} }} \leq
n^{-{n^{\frac{1}{2} - o\left( 1 \right)}}}.
\end{eqnarray*}

By  the Union Bounds, with probability at least $1-\frac{\log
n}{101}o(n^{ - {n^{\frac{1}{2} - o\left( 1 \right)}}}) \geq
1-o(n^{-n^{2/5}})$, we have that for every $T_{v_j}$, the number of
edges between $L(T_{u_i})$  and $L(T_{v_j})$ is at most
$8\mathbb{E}\left[ {{A_{ij}}} \right]$=$8p{\left( {\frac{{\log
n}}{{101}}} \right)^{2\left( {\frac{3}{4} - \varepsilon } \right)D -
2}}$.

Therefore, with probability at least $1-o(n^{ - {n^{\frac{1}{2} -
o\left( 1 \right)}}})-o(n^{-n^{2/5}})=1-o(n^{-n^{2/5}})$ there are
at least $\frac{\frac{100}{101}\left( {\frac{{\log n}}{{101}}}
\right)^{2\left( {\frac{3}{4} - \varepsilon } \right)D - 1} \cdot
p}{8{\left( {\frac{{\log n}}{{101}}} \right)^{2\left( {\frac{3}{4} -
\varepsilon } \right)D - 2}} \cdot p}$$=\frac{100}{101^2\times
8}\log n$ vice trees $T_{v_j}$, such that vertices in $L(T_{u_i})$
and $L(T_{v_j})$ can be connected with edges. Moreover, using the
Union Bounds, with probability at least $1-\frac{\log
n}{101}o(n^{-n^{2/5}})$ $\geq 1-o(n^{-n^{1/5}})$, each vice tree of
$T_{u}$ can be connected to $\frac{100}{101^2\times 8}\log n$ vice
trees of $T_v$ with edges. Hence there are at least
$\frac{100}{101^2\times 8}\log n$ pairs $\{T_{u_i},T_{v_j}\}$ such
that vertices of $L(T_{u_i})$ and $L(T_{v_j})$  can be connected by
edges.

For convenience, let $a\log n= \frac{100}{101^2\times 8}\log n$.
Without loss of generality, assume these $a\log n$ pairs be
$T_{u_\ell}$ and $T_{v_\ell}$, $\ell=1,2,\cdots, a\log n$. Now we
show that, for the remaining large vertex $w$, we can find at least
$k$ internally disjoint trees connecting $u$, $v$ and $w$.

Note that we can assume that $w$ is not in $T_u$ and $T_v$, since
otherwise we can prune the tree by deleting the subtree rooted at
$w$ (just like the way to deal with small vertices), and we can
still obtain $\frac{\log n}{101}$-ary trees rooted at $u$ and $v$,
respectively.

With the similar argument in the proof of Lemma \ref{lem6}, we can
construct a $\frac{\log n}{101}$-ary tree $T_w$ of depth $(\frac{1}{4} + 2\varepsilon)D $
rooted at $w$, and $T_u$, $T_v$, $T_w$ are vertex disjoint. Note that at this time the number of
small vertices in each level is at most one.

Let $w_1,\ldots, w_{\frac{\log n}{101}}$ be the vertices of the
first depth of $T_w$. Let $Q_i=T_{u_i}\cup T_{v_i}$, $i=1,2,\ldots ,
a\log n$. Then $|Q_i|>2{\left( {\frac{{\log n}}{{101}}}
\right)^{\left( {\frac{3}{4} - \varepsilon } \right)D-1 }}$. For any
fixed $Q_j$, let $q_j$ denote the probability that there~exists~at
least one edge between $T_{w_j}$ and $Q_j$. Then
\begin{eqnarray*}
q_j &=&   1-\Pr[\mathrm{there~is~no~edge~between}~T_{w_j}~\mathrm{and}~Q_j]\\
\null  &=& 1 - {\left( {1 - p} \right)^{|{Q_j}| \cdot |{T_{{w_j}}}|}}\\
\null  &\geq& 1 - {\left( {1 - p} \right)^{2{{\left( {\frac{{\log n}}{{101}}}
\right)}^{\left( {\frac{3}{4} - \varepsilon } \right)D - 1}} \cdot {{
\left( {\frac{{\log n}}{{101}}} \right)}^{\left( {\frac{1}{4} + 2\varepsilon } \right)D - 1}}}}\\
\null &\geq& 1 - {e^{ - 2\frac{{\log n}}{n}{{\left( {\frac{{\log n}}{{101}}}
\right)}^{\left( {1 + \varepsilon } \right)D -2}}}}\\
\null  &=& 1 - {e^{ - 2\log n \cdot \frac{{{{101}^2}}}{{{{\log }^2}n}} \cdot{n^{\varepsilon
- \frac{{4.6\left( {1 + \varepsilon } \right)}}{{\log \log n}}}} }}\\
\null  &=& 1 - {e^{ - 2\frac{{{{101}^2}}}{{\log n}}\cdot
{n^{\varepsilon  - \frac{{4.6\left( {1 + \varepsilon }
\right)}}{{\log \log n}}}}  }}.
\end{eqnarray*}

Since for  any $i\neq j$, $q_i=q_j$, let $q=q_i=q_j$, and let
$\mathcal{A}$ be the event that there are at most $k-1$ pairs of
$\{T_{w_\ell},Q_\ell\}$, such that there exist edges between
$T_{w_\ell}$ and $Q_\ell$, where $\ell=1,\ldots, a\log n$.

Consider the upper bound of the probability that  $\mathcal{A}$
happens, we can deduce that
\begin{eqnarray*}
\Pr[\mathcal{A}]  &\leq &\sum\limits_{i = 0}^{k - 1} {\left(
{\begin{array}{*{20}{c}}
{a\log n}\\
i
\end{array}} \right)} {q^i}{\left( {1 - q} \right)^{a\log n - i}}\\
\null &\leq& k{\left( {\frac{{a\log n \cdot e}}{{k - 1}}}
\right)^{k - 1}}{\left( {1 - q} \right)^{\frac{a}{2}\log n }}\\
\null &<& k{\left( {\frac{{a\log n \cdot e}}{{k - 1}}}
\right)^{k - 1}}{e^{ - 2\frac{{{{101}^2}}}{{\log n}}
\cdot {n^{\varepsilon  - \frac{{4.6\left( {1 + \varepsilon }
\right)}}{{\log \log n}}}} \cdot \frac{a}{2}\log n}}\\
\null &=& k{\left( {\frac{{a\log n \cdot e}}{{k - 1}}}
\right)^{k -
1}}{n^{ - \frac{{100}}{8}\frac{{{n^{^{\varepsilon  -
\frac{{4.6\left( {1 + \varepsilon } \right)}}{{\log \log
n}}}}}}}{{\log n}}}}< n^{-10}.
\end{eqnarray*}
This indicates  there are at least $k$ pairs of
$\{T_{w_\ell},Q_\ell\}$, such that $T_{w_\ell}$ and $Q_\ell$ can be
connected by edges, a.s..

Without loss of generality, assume that there exists  edges
connecting $T_{w_\ell}$ and $Q_\ell$, where $\ell=1,\cdots, k$. Now
we will construct $k$ internally disjoint trees connecting $u$, $v$
and $w$. For each $i$ with $i=1,\cdots, k$, suppose that $w'\in
V(T_{w_i})$ is adjacent to $x\in V(T_{u_i})$, edge $yz$ connects
$L(T_{u_i})$ and $L(T_{v_i})$, where $y \in L(T_{u_i})$, $z \in
L(T_{v_i})$. Let $P_{T_{u_i}}(x,y)$ denote the path connecting $x$
and $y$ in $T_{u_i}$.

$\bullet$ If $u_i$ is  not in $P_{T_{u_i}}(x,y)$, we construct the
tree $\{uu_i\}\cup P_{T_{u_i}}(u_i,x)\cup P_{T_{u_i}}(x,y)\cup
\{yz\}\cup P_{T_{v_i}}(z,v_i)\cup \{v_iv\}\cup
P_{T_{w_i}}(w',w_i)\cup \{w_iw\}$.

$\bullet$ If $u_i$ is  contained in $P_{T_{u_i}}(x,y)$, we construct
the tree $\{uu_i\}\cup  P_{T_{u_i}}(x,y)\cup \{yz\}\cup
P_{T_{v_i}}(z,v_i)\cup \{v_iv\}\cup P_{T_{w_i}}(w',w_i)\cup
\{w_iw\}$.

For the case that $x\in V(T_{v_i})$, the tree connecting $u$, $v$
and $w$ can be constructed similarly.

Thus we construct $k$ trees  connecting  large vertices $u$, $v$ and
$w$, and it is easy to get that all these trees are internally
disjoint.

Now we deal with small vertices. From the previous argument, if the
three given vertices are all large, we can find at least $k$
internally disjoint trees connecting them. So we assume that there
are at least one small vertex of the given vertices $u$, $v$ and
$w$. By Lemma 4, it is easy to obtain the following three facts.

\begin{enumerate}
  \item The neighbors of a small vertex are large vertices, a.s..
  \item A large vertex can have at most one neighbor that is small, a.s..
  \item Any two small vertices have no common neighbors, a.s..
\end{enumerate}

Combining the three facts above, we can take $k$ large neighbors of
$u$, $v$, $w$, denoted by $u_1, \ldots, u_{k}$, $v_1, \ldots, v_{k}$
and  $w_1, \ldots, w_{k}$, respectively, and all these $3k$ vertices
are different.

Firstly, using the method described before, we can find a tree $T_1^*$
connecting $u_1$, $v_1$ and $w_1$. Note that the number of edges in
$T_1^*$ is at most $3\left( {\frac{3}{4} - \varepsilon } \right)D +
\left( {\frac{1}{4} + 2\varepsilon } \right)D = \left( {\frac{5}{2}
- \varepsilon } \right)D$.

Then we find a tree $T_2^*$ connecting $u_2$, $v_2$ and $w_2$. In
order to make $T_1^*$ and $T_2^*$ internally disjoint,  when we
construct BFS tree rooted at $u_2,~v_2$ and $w_2$, we treat the
edges with one endpoint in $V(T_1^*)$ as bad edges, too. By Remark
\ref{rem2} and the similar argument as we deal with BFS tree rooted
at a large vertex, we can find a tree $T_2^*$ connecting $u_2$,
$v_2$ and $w_2$.

Continue that process, until we find trees $T_j^*$ connecting $u_j$,
$v_j$ and $w_j$ for all $j=1,\cdots, k$.

Let $T_i= T_i^* \cup \{uu_i, vv_i, ww_i\} $ for $i=1,\cdots, k$.
Apparently, these are $k$ internally disjoint connecting $u$, $v$
and $w$. Thus, for any fixed three vertices $u$, $v$, $w$, with
probability at least
$1-o(n^{-3})-o(n^{-n^{1/5}})-o(n^{-10})=1-o(n^{-3})$, we can find
$k$ vertex disjoint connecting them.

Consequently, for all possible three vertices $u$, $v$, $w$, by the
Union Bounds, we can find $k$ internally disjoint connecting them
with probability at least $1-n^3\cdot o(n^{-3})=1-o(1)$.\qed

In conclusion, combining the results of Theorems \ref{thm2} and
\ref{thm3}, we can derive Theorem \ref{thm1} immediately.

\vspace{3ex} Considering the generalized edge-connectivity of random
graphs, it is easy to get that $\kappa_3(G) \leq \lambda_3(G)\leq
\delta (G)$, for any connected graph $G$. Furthermore, Li et al.
\cite{lm7} gave the result as follows.
\begin{lem}[\cite{lm7}]\label{lem8}
Let $G$ be a graph  of order $n$, then $\lambda_r(G)\leq \lambda(G)$,
for any $3\leq r\leq n$. Moreover, the upper bound is sharp.\qed
\end{lem}

Combining Theorems \ref{lem9} and \ref{thm3}, Lemmas \ref{lem7} and
\ref{lem8}, we can get the result of Corollary \ref{cor1}.

\end{document}